\documentstyle[11pt,fullpage,twoside,amssym,saks]{article}

\def\XXint#1#2#3{{\setbox0=\hbox{$#1{#2#3}{\int}$}
     \vcenter{\hbox{$#2#3$}}\kern-.5\wd0}}

\begin{document}
\markboth{\centerline{R. SALIMOV AND E. SEVOST'YANOV}}
{\centerline{$ACL$\&DIFFERENTIABILITY...}}

\def\kohta #1 #2\par{\par\noindent\rlap{#1)}\hskip30pt
\hangindent30pt #2\par}
\def\esssup{\operatornamewithlimits{ess\,sup}}
\def\tomes{\mathop{\longrightarrow}\limits^{mes}}
\def\ts{\textstyle}
\def\I{\roman{Im}}
\def\mes{\mbox{\rm mes}}
\def\Rm{{{\Bbb R}^m}}
\def\Rn{{{\Bbb R}^n}}
\def\Rk{{{\Bbb R}^k}}
\def\R3{{{\Bbb R}^3}}
\def\lR{{\overline {{\Bbb R}}}}
\def\lRn{{\overline {{\Bbb R}^n}}}
\def\lRm{{\overline {{\Bbb R}^m}}}
\def\lBn{{\overline {{\Bbb B}^n}}}
\def\Bn{{{\Bbb B}^n}}
\def\R{{\Bbb R}}
\def\Z{{\Bbb Z}}
\def\C{{\Bbb C}}
\def\B{{\Bbb B}}
\def\e{{\varepsilon}}
\def\L{{\Lambda}}
\def\f{{\varphi}}
\def\t{{\tau}}
\def\x{{\chi}}
\def\d{{\delta }}
\def\b{{\beta }}
\def\D{{\Delta }}
\def\c{{\circ }}
\def\tg{{\tilde{\gamma}}}
\def\a{{\alpha }}
\def\p{{\psi }}
\def\P{{\Psi }}
\def\F{{\frak{F}}}
\def\m{{\mu }}
\def\r{{\rho }}
\def\O{{\Omega }}
\def\s{{\sigma }}
\def\l{{l}\, }
\def\g{{\gamma }}
\def\G{{\Gamma }}
\def\D{{\Delta }}
\let\text=\mbox
\let\Cal=\cal

\def\cc{\setcounter{equation}{0}
\setcounter{figure}{0}\setcounter{table}{0}}

\overfullrule=0pt

\def\eqb{\begin{equation}}
\def\eqe{\end{equation}}
\def\eb{\begin{eqnarray}}
\def\ee{\end{eqnarray}}
\def\ebnn{\begin{eqnarray*}}
\def\eenn{\end{eqnarray*}}
\def\db{\begin{displaystyle}}
\def\de{\end{displaystyle}}
\def\tb{\begin{textstyle}}
\def\te{\end{textstyle}}
\def\exb{\begin{ex}}
\def\exe{\end{ex}}
\def\bth{\begin{theo}}
\def\eth{\end{theo}}
\def\bcor{\begin{corol}}
\def\ecor{\end{corol}}
\def\blem{\begin{lemma}}
\def\elem{\end{lemma}}
\def\brem{\begin{rem}}
\def\erem{\end{rem}}
\def\bpr{\begin{propo}}
\def\epr{\end{propo}}
\title{{\bf $ACL$  AND DIFFERENTIABILITY OF THE OPEN DISCRETE RING $(p,Q)$-MAPPINGS}}

\author{{\bf R. Salimov and E. Sevost'yanov}\\}
\date{\today \hskip 4mm }
\maketitle

\large \begin{abstract} We study the so--called ring $Q$--mappings
which are the natural generalization of quasiregular mappings. It is
proved that open discrete ring $Q$--mappings are differentiable a.e.
and belong to the class $ACL$ in ${\Bbb R}^n$, $n\ge 2$,
furthermore, $f\in W_{loc}^{1,1}$ provided that $Q\in L^{1}_{loc}$.
\end{abstract}


\section{Introduction}\label{s1}

Recall that, given a family of paths $\Gamma$ in $\Rn$, a Borel
function $\varrho:\Rn\to[0,\infty]$ is called {\bf admissible} for
$\Gamma$, abbr. $\varrho\in adm\,\Gamma$, if \eqb\label{eq13.2}
\int\limits_{\gamma}\varrho\,ds\ \geq\ 1\eqe for all
$\gamma\in\Gamma$. The  {\bf modulus} of $\Gamma$ is the quantity
\eqb\label{eq13.3}M_{\alpha}(\Gamma)\ =\ \inf\limits_{\varrho\in
adm\,\Gamma}\int\limits_{G}\varrho^{\alpha}(x)\ dm(x)\ .\eqe\
\medskip

Let $D$ be a domain in ${\Bbb R}^n,\,\,n\ge 2,$ and
$f:D\rightarrow {\Bbb R}^n$  be a $Q$--qua\-si\-con\-for\-mal
mapping. Then necessarily
\begin{equation}\label{eq1}
M_n\left(f\Gamma\right)\le\int\limits_{D} K_I(x,f)\cdot\rho^n(x)\
dm(x)
\end{equation}
for every family $\Gamma$ of paths in $D$ and every admissible
function $\rho$ for $\Gamma,$ see e.g. \cite{BGMV}, where $K_I(x,f)$
stands for the well-known inner or outer dilatation of $f$ at $x.$
One can replace the above necessary condition  with the following,
equivalent by Gehring's result \cite{Ge}, inequality
\begin{equation}\label{eq4}
M_n\left(f\left(\Gamma\left(S_1,\,S_2,\,A\right)\right)\right)\le\int\limits_{A(r_1,
r_2, x_0)} K_I(x,f)\cdot\rho^n\left(|x-x_0|\right)\ dm(x)
\end{equation}
for every point $x_0\in D$  and every $r_1, r_2,$ such that
$0<r_1<r_2<r_0\,=\,{\rm dist\,} \left(x_0,\partial D\right),$ where
$A=A(r_1, r_2, x_0) = \{ x\,\in\,{\Bbb R}^n : r_1<|x-x_0|<r_2\}$ and
$\Gamma\left(S_1,\,S_2,\,A\right)$ is a family of all paths joining
the spheres $ S_{\,i}\,=\,S(x_0,r_i) = \{ x\,\in\,{\Bbb R}^n :
|x-x_0|=r_i\},$ $i=1,2,$ in $A(r_1, r_2, x_0).$ The above
inequalities  together with the modulus technique are the powerful
tools for the study of quasiconformal (quasiregular) mappings in the
plane and in space, see e.g. \cite{GL}, \cite{Va}, \cite{Re$_2$} and
\cite{Ri}. In order to extend as much as possible the set of maps
for the study of which the well developed  modulus technique can be
also applied, we replace in (\ref{eq1}) (in (\ref{eq4})) the
dilatation $K_I(x,f)$ with a measurable function $Q(x),$ say of the
class $L_{loc}^1(D),$ and then declare the inequality
\begin{equation}\label{eq13.1}
M_p(f\Gamma)\leq\int\limits_{D} Q(x)\cdot\rho^p(x)\ dm(x)
\end{equation}
or
\begin{equation}\label{eq5}
M_p\left(f\left(\Gamma\left(S_1,\,S_2,\,A\right)\right)\right)\
\leq \int\limits_{A} Q(x)\cdot \eta^p(|x-x_0|)\ dm(x)\,,
\end{equation} due to
\cite{MRSY$_2$}, as the necessary condition for the mapping
$f:D\rightarrow {\Bbb R}^n$  to belong to the class of the
$Q$--homeomorphisms,  or the ring $Q$--homeomorphisms, respectively,
etc. See also the conception of the weighted modulus,
\cite{AC$_1$}--\cite{AC$_2$}, and applications of the
$Q$--homeomorphisms, cf. \cite{Cr$_1$}--\cite{Cr$_2$}, \cite{BGR}
and \cite{S$_1$}--\cite{S$_2$}.

\bigskip
Note that, if $f$ is homeomorphism, the inequality (\ref{eq5}) holds
at every point $x_0\in D$ and $Q(x)\le K$ a.e., the definitions of
$Q$--homeomorphism and ring $Q$--homeomorphism are equivalent, see
\cite{Ge}, and give that $f$ is $K$--qua\-si\-con\-for\-mal mapping.
Moreover, every $K$--qua\-si\-con\-for\-mal (or $K$--quasiregular)
mapping satisfies to (\ref{eq5}) and (\ref{eq13.1}) with $Q(x)\equiv
K.$ Our paper is devoted to the study of mappings having unbounded
$Q(x)$ in above definitions.

\bigskip
Recall that a mapping $f:D\rightarrow {\Bbb R}^n$ is said to be {\bf
absolutely continuous on lines}, write $f \in ACL, $ if all
coordinate functions $f = \left(f_1,\ldots,f_n\right)$ are
absolutely continuous on almost all straight lines parallel to the
coordinate axes for any $n$--dimensional parallelepiped $P$ with
edges parallel to the coordinate axes and such that
$\overline{P}\subset D.$

\bigskip
It is well--known that quasiconformal and quasiregular mappings are
absolutely continuous on lines, see e.g. Corollary 31.4 in
\cite{Va}, Lemma 4.11 and Theorem 4.13 in \cite{MRV$_1$}, and
differentiable a.e., see e.g. Corollary 32.2 in \cite{Va}, Theorem
2.1 Ch. I in \cite{Ri} and Theorem 4 in \cite{Re$_1$}. Moreover, in
the plane case, every $ACL$--homeomorphism is differentiable a.e.,
see \cite{GL}. However, above results did not give any information
about differentiability (or $ACL$) for more general mappings having
non--bounded dilatation. The first steps in this direction were made
in the work of one of authors, see \cite{Sal}. More detail, it has
been shown that $Q$--homeomorphisms are differentiable a.e. and
belong to the class $ACL$ provided that a function $Q$ is locally
integrable. In the present paper we extend these results to open
discrete mappings satisfying the conditions of the type (\ref{eq5}).

\bigskip
Thus, the goal of the present paper is to prove the following:

\bigskip
I. Open discrete ring $(p,Q)$--mappings $f:D\rightarrow
\overline{{\Bbb R}^n}$ with $Q\in L^1_{loc}$ and $p>n-1$ are
differentiable a.e. in $D.$

II. Open discrete ring $(p,Q)$--mappings $f:D\rightarrow
\overline{{\Bbb R}^n}$ with $Q\in L^1_{loc}$ and $p>n-1$ belongs
to the class $ACL$ in $D.$

III. Open discrete ring $(p,Q)$--mappings $f:D\rightarrow
\overline{{\Bbb R}^n}$ with $Q\in L^1_{loc}$ and $p>n-1$  belong
to the Sobolev class $W_{loc}^{1,1}$ and satisfy to the inequality
$$\Vert f^{\,\prime}(x)\Vert^p\ \le\ C_\cdot|J(x,f)|^{1-n+p}\ Q^{n-1}(x)$$
a.e. where a constant $C$ depends only on $n$ and $p.$

\section{Preliminaries}\cc

Let $D$ be a domain in ${\Bbb R}^n,\,\,n\ge 2$. A mapping
$f:D\rightarrow {\Bbb R}^n$ is said to be {\bf discrete} if the
preimage $f^{-1}\left(y\right)$ of every point $y\,\in\,{\Bbb R}^n$
consists of isolated points, and an {\bf open} if the image of every
open set $U\,\subseteq\,D$ is open in ${\Bbb R}^n\,.$ The notation
$G\Subset D$ means that $\overline{G}$ is a compact subset of $D.$
We suppose that $f:D\,\rightarrow\,{\Bbb R}^n$ is continuous and
sense--preserving, i.e. a topological index $\mu\,
(y,\,f,\,G)\,>\,0$ for any $G\Subset D$ and $y\,\in\,f(G)\setminus
f\left(\partial G\right).$  A {\bf neighborhood} of a point $x$ or a
set $A$ is an open set containing $x$ or $A,$ correspondingly.
Suppose that $x\in D$ has a connected neighborhood $G$ such that
$\overline{G}\cap f^{\,-1}\left(f(x)\right)=\left\{x\right\}.$ Then
$\mu\, \left(f(x),\,f,\,G\right)$ is well--defined and independent
of the choice of $G$ for discrete open $f$ and denoted by $i(x,f).$
For $f:D\,\rightarrow\,{\Bbb R}^n$ and $E\subset D,$ we use the {\bf
multiplicity functions}
$$N(y,f,E)\,=\,{\rm card}\,\left\{x\in E: f(x)=y\right\}\,,$$
$$N(f,E)\,=\,\sup\limits_{y\in{\Bbb R}^n}\,N(y,f,E)\,.$$
In what follows, we also use the notations
$B(x_0,\,r)\,=\,\left\{x\in{\Bbb R}^n: \vert
x-x_0\vert\,<\,r\right\}$ and $\overline{{\Bbb R}^n}={\Bbb
R}^n\cup\{\infty\}.$ The  above definitions can be extended in a
natural way to mappings $f:D\rightarrow \overline{{\Bbb R}^n}$.

The following notion is motivated by the Gehring ring definition
of quasiconformality, see \cite{Ge}, and generalizes a notion of
ring $Q$--homeomorphism, see \cite{RSY}.

Given a domain $D$ and two sets $E$ and $F$ in ${\overline{{\Bbb
R}^n}},$ \,\,$n\ge 2,$\,\, $\G (E,F,D)$ denotes the family of all
paths $\g:[a,b] \rightarrow {\overline{{\Bbb R}^n}}$ which join $E$
and $F$ in $D$, i.e., $\g(a) \in E, \ \g(b) \in F$ and $\g(t) \in D$
for $a<t<b$. We set $\G(E,F)= \G(E,F,{\overline{{\Bbb R}^n}})$ if
$D={\overline{{\Bbb R}^n}}.$ Let $r_0\,=\,{\rm dist\,}
(x_0\,,\partial D)$ and $Q:D\rightarrow\,[0\,,\infty]$ is a
measurable function. Set
$$ A(r_1,r_2,x_0) = \{ x\,\in\,{\Bbb R}^n : r_1<|x-x_0|<r_2\}\,, $$
$$ S_{\,i}\,=\,S(x_0,r_i) = \{ x\,\in\,{\Bbb R}^n :
|x-x_0|=r_i\}\,\,,\ \ \ i=1,2. $$
A homeomorphism $f:D\rightarrow \overline{{\Bbb R}^n}$ is said to
be a {\bf ring  $(\p,Q)$--homeomorphism at a point $x_0\,\in\,D$,}
if
\begin{equation}\label{eq2} M_{p}\left(f\left(\Gamma\left(S_1,\,S_2,\,A\right)\right)\right)\ \leq
\int\limits_{A} Q(x)\cdot \eta^{p}(|x-x_0|)\ dm(x)
\end{equation}
holds for every annulus $A=A(r_1,r_2, x_0),$\, $0<r_1<r_2< r_0$ and
every measurable function $\eta : (r_1,r_2)\rightarrow [0,\infty
]\,$ such that
$$
\int\limits_{r_1}^{r_2}\eta(r)\ dr\ \ge\ 1\,.
$$
If $(\ref{eq2})$ holds for every $x_0\,\in\,D\,,$  $f$ is said to
be a ring  $(p,Q)$--homeomorphism. In general case, every
 $(p,Q)$--homeomorphism $f:D\rightarrow \overline{{\Bbb R}^n}$ is a ring
 $(p,Q)$--homeomorphism, but the inverse conclusion, generally speaking,
is not true. In \cite{RSY} there are examples of ring
Q-homeomorphisms in a fixed point $x_0$ such that $Q(x)\in(0, 1)$
on some set for which $x_0$ is a density point. We will not
discuss here these connections in more details.

Let $D\subset{{\Bbb R}^n}\,,$\,\,$n\ge 2\,,$ be a domain and
$Q:D\rightarrow\,[0\,,\infty]$ be a measurable function. We say
that a continuous sense--preserving mapping $f:D \rightarrow
\overline{{\Bbb R}^n}$ is a {\bf ring  $(p,Q)$mapping in $D$} if
$(\ref{eq2})$ holds for every $x_0\,\in\,D.$ Note that
correspondingly to these definitions the class of the so--called
 $(p,Q)$--mappings which consists of the continuous sense--preserving
mappings satisfying the condition (\ref{eq13.1}) is included in
the class of ring  $(p,Q)$--mappings. Thus, all results for ring
$(p,Q)$--mappings formulated below hold, in particular, for
$(p,Q)$--mappings.

Correspondingly to \cite{MRV$_1$} a {\bf condenser} is a pair
$E\,=\,(A, C)$ where $A\subset {\Bbb R}^n$ is open and $C$ is
non--empty compact set contained in A\,. A condenser $E\,=\,(A,
C)$ is said to be in a domain $G$ if $A\subset G\,.$ For a given
condenser $E\,=\,\left(A,\,C\right),$ we set

\begin{equation}\label{equ5}
{\rm cap_{p}}\,E\,\,=\,\,{\rm
cap_p}\,\left(A,\,C\right)\,\,=\,\,\inf
\limits_{u\,\in\,W_0\left(E\right)
}\,\,\int\limits_{A}\,\vert\nabla u\vert^{p}\,\,dm(x)
\end{equation}
where $W_0(E)\,=\,W_0(A,\,C)$ is the family of non--negative
functions $u:A\rightarrow R^1$ such that (1)\,\,$u$ is continuous
and finite on $A,$\,\,(2)\,\, $u(x)\ge 1$ for $x\in C,$ and
(3)\,\,$u\,$ is $\,ACL.$ In the above formula
$$\vert\nabla
u\vert\,=\,{\left(\sum\limits_{i=1}^n\,{\left(\partial_i u\right)}^2
\right)}^{1/2}\,.$$
The quantity ${\rm cap_p\,}E$ is called the {\bf p-capacity} of
the condenser $E.$

We say that a family of curves $\Gamma_1$ is minorized by a family
$\Gamma_2,$  denoted by $\Gamma_1>\Gamma_2,$ if for every curve
$\gamma\in \Gamma_1$ there is a subcurve that belongs to the
family $\Gamma_2.$ It is known that $M_{p}(\Gamma_1)\le
M_{p}(\Gamma_2)$ as $\Gamma_1>\Gamma_2,$ see Theorem 6.4 in
\cite{Va}.

\section{Differentiability}\cc

Let $f:D \rightarrow {\Bbb R}^n$ be a discrete open mapping. Let
$\beta: [a,\,b)\rightarrow {\Bbb R}^n$ be a path and
$x\in\,f^{-1}\left(\beta(a)\right).$ A path $\alpha:
[a,\,c)\rightarrow D$ is called a {\bf maximal $f$--lifting} of
$\beta$ starting at $x$ if $(1)\,\,\,\, \alpha(a)\,=\,x\,;$
$(2)\,\,\,\, f\circ\alpha\,=\,\beta|_{[a,\,c)}\,;$ $(3)$\,\,\,if
$c\,<\,c^{\prime}\,\leq\,b,$ then there is no path $\alpha^{\prime}:
[a,\,c^{\prime})\rightarrow D$ such that
$\alpha\,=\,\alpha^{\prime}|_{[a,\,c)}$ and $f\circ
\alpha^{\,\prime}\,=\,\beta|_{[a,\,c^{\prime})}.$ If $f$ is a
discrete open mapping, then every path $\beta$ with
$x\in\,f^{-1}\left(\beta(a)\right)$ has  a maximal $f$--lifting
starting at a point $x,$ see Corollary $3.3$ Ch.II in \cite{Ri}. We
need the following statement, see Proposition $10.2$ Ch. II in
\cite{Ri}.
\begin{lemma}{}\label{lem2.2}
Let $E\,=\,(A,\,C)$ be a condenser in ${\Bbb R}^n$ and let
$\Gamma_E$ be the family of all paths of the form
$\gamma:\,[a,\,b)\,\rightarrow\,A$ with $\gamma(a)\,\in\,C$ and
$\vert\gamma\vert\cap\left(A\setminus F\right)\,\neq\,\varnothing$
for every compact $F\,\subset\,A.$ Then ${\rm
cap_p}\,E\,=\,M_p\left(\Gamma_E\right).$
\end{lemma}

\begin{theo}{}\label{th1} Let $D$ be a domain in ${\Bbb R}^n$, $n\ge 2$,
and $f:D\rightarrow \overline{{\Bbb R}^n}$ be a ring
$(p,Q)$--mapping with $Q\in L^1_{loc}$ and $p>n-1$. Suppose that
$f$ is discrete and open. Then $f$ is differentiable a.e. in $D$.
\end{theo}

{\it Proof.} Without loss of generality we may assume that
$\infty\notin D^{\,\prime}=f(D).$ Let us consider the set function
$\Phi(B)=m\left(f(B)\right)$ defined over the algebra of all the
Borel sets $B$ in $D.$ By 2.2, 2.3 and 2.12 in \cite{MRV$_1$}
\begin{equation}\label{3.1}
\varphi(x)=\limsup\limits_{\varepsilon\rightarrow
0}\frac{\Phi(B(x,\varepsilon))}{\Omega_{n}\varepsilon^{n}}<\infty
\end{equation}
for a.e. $x\in D.$ Consider the spherical ring
$R_\varepsilon(x)=\{y:\ \varepsilon<|x-y|<2\varepsilon\}$, $x\in
D,$ with $\varepsilon>0$  such that $B(x,2\varepsilon)\subset D$.
Note that
$E=\left(B\left(x,2\varepsilon\right),\overline{B\left(x,\varepsilon\right)}\right)$
is a condenser in $D$ and
$fE\,=\,\left(fB\left(x,2\varepsilon\right),f\overline{B\left(x,\varepsilon\right)}\right)$
is a condenser in $D^{\,\prime}.$ Let $\Gamma_E$ and $\Gamma_{fE}$
be path families from Lemma \ref{lem2.2}. Then
\begin{equation}\label{eq5aaa}
{\rm cap_p}\ \left(fB(x,2\varepsilon\right), f\overline{B(x
,\varepsilon)})=M_p\left(\Gamma_{fE}\right).
\end{equation}
Let $\Gamma^{*}$ be a family of maximal $f$--liftings of
$\Gamma_{fE}$ starting at $\overline{B\left(x,\varepsilon\right)}.$
We show that $\Gamma^{*}\,\subset\,\Gamma_E\,.$ Suppose the
contrary. Then there is a path $\beta: [a,\,b)\rightarrow {\Bbb
R}^n$ of $\Gamma_{fE}$ such that the corresponding maximal
$f$--lifting $\alpha:[a,\,c)\rightarrow B(x,2\varepsilon)$ is
contained in some compact $K$ inside of $B(x,2\varepsilon)\,.$ Thus
$\overline{\alpha}$ is a compactum in $B(x,2\varepsilon),$ see
Theorem 2, \,$\S\, 45$ in \cite{Ku}. Remark that $c\neq\,b.$ Indeed,
in the contrary case $\overline{\beta}$ is a compact in $f(A)$ that
contradicts to the condition $\beta\,\in\,\Gamma_{fE}.$ Consider the
set
$$G\,=\,\left\{x\in {\Bbb R}^n: x\,=\,\lim\limits_{k\rightarrow\,\infty} \alpha(t_k)
 \right\}\,,\,\,\,\,\,\,\,\,\,t_k\,\in\,[a,\,c)\,,\,\,
 \lim\limits_{k\rightarrow\infty}t_k\,=\,c\,.$$
Without loss of generality we may assume that $t_k$ is the monotone
sequence. By continuity of $f,$ for $x\in\,G\,,$
$f\left(\alpha(t_k)\right)\rightarrow\,f(x)$ as
$k\rightarrow\,\infty$ where $t_k\,\in\,[a,\,c),\,t_k\rightarrow c$
as $k\rightarrow \infty\,.$ However,
$f\left(\alpha(t_k)\right)\,=\,\beta(t_k)\,\rightarrow\,\beta(c)$ as
$k\rightarrow\infty\,.$ Thus, $f$ is a constant in $G\subset
B(x,2\varepsilon)\,.$ On the other hand, from the Cantor condition
on the compact $\overline{\alpha},$
%
%
$$G\,=\,\bigcap\limits_{k\,=\,1}^{\infty}\,\overline{\alpha\left(\left[t_k,\,c\right)\right)}\,=\,
\limsup\limits_{k\rightarrow\infty}\alpha\left(\left[t_k,\,c\right)\right)\,=\,
\liminf\limits_{k\rightarrow\infty}\alpha\left(\left[t_k,\,c\right)\right)\,\neq\,\varnothing
$$
%
by monotonicity of the sequences of connected sets
$\alpha\left(\left[t_k,\,c\right)\right),$ see \cite{Ku}. Thus,
$G$ is connected by $I (9.12)$ in \cite{Wh}. Consequently, $G$ is
a single point by discreteness of $f.$ So a path
$\alpha:[a,\,c)\rightarrow\,B(x,2\varepsilon)$ can be extended to
$\alpha:[a,\,c]\rightarrow K\,\subset\,B(x,2\varepsilon)$ and
$f\left(\alpha(c)\right)\,=\,\beta(c)\,.$ By Corollary $3.3$ Ch.
$II$ in \cite{Ri} we can construct a maximal $f$--lifting
$\alpha^{\,\prime}$ of $\beta|_{[c,\,b)}$ started at $\alpha(c).$
United the liftings $\alpha$ and $\alpha^{\,\prime},$ we have a
new $f$--lifting $\alpha^{\,\prime\prime}$ of $\beta$ defined on
$[a, c^{\prime}),$ \,\,$c^{\,\prime}\,\in\,(c,\,b),$ that
contradicts to the maximality of $f$--lifting $\alpha.$ Thus
$\Gamma^{*}\,\subset\,\Gamma_E\,.$ Remark that
$\Gamma_{fE}\,>\,f\Gamma^{*}$  and, consequently,
%
$$M_p\left(\Gamma_{fE}\right)\le M_p\left(f\Gamma^{*}\right)\le
M_p\left(f\Gamma_E\right)\,.$$

Let $\left\{r_i\right\}_{i=1}^{\infty}$ be an arbitrary sequence
of numbers with $\varepsilon<r_i<2\varepsilon$ such that
$r_i\rightarrow 2\varepsilon-0.$ Denote by $\Gamma_i$ a family of
paths joining the spheres $|x|=\varepsilon$ and $|x|=r_i$ in a
ring $\varepsilon<|x|<r_i.$ Then $\Gamma_E>\Gamma_i$ for every
$i\in {\Bbb N}.$
Consider the family of functions
$$\eta_{i,\varepsilon}(t)\,=\,\left \{\begin{array}{rr}
\frac{1}{r_i-\varepsilon}, & {\rm если   } \ t\in (\varepsilon,
r_i),
\\
0, & {\rm если } \ t\in {\Bbb R}\setminus (\varepsilon, r_i)\,.
\end{array}\right.
$$
By definition of a ring $Q$--mapping
\begin{equation}\label{eq9}
M_p(f\Gamma_E)\le M_p(f\Gamma_i)\le
\frac{1}{(r_i-\varepsilon)^p}\int\limits_{\varepsilon<|x|<r_i}Q(x)\,dm(x)\le
\frac{1}{(r_i-\varepsilon)^p}\int\limits_{B(x,2\varepsilon)}Q(x)\,dm(x)\,.
\end{equation}
Letting to the limit in (\ref{eq9}) as $i\rightarrow\infty$, we
obtain
\begin{equation}\label{eq10}
M_p(f\Gamma_E)\le
\frac{1}{\varepsilon^p}\int\limits_{B(x,2\varepsilon)}Q(x)\,dm(x)\,.
\end{equation}
From (\ref{eq5aaa}) and (\ref{eq10})
\begin{equation}\label{2.4}{\rm cap_p}\ (fB(x,2\varepsilon),f\overline{B(x,\varepsilon)})\leq\,
\frac{1}{\varepsilon^p}\int\limits_{B(x,2\varepsilon)}Q(x)\,dm(x)\,.
\end{equation}
On the other hand, by Proposition 6 in \cite{Kr}
\begin{equation}\label{2.5} {\rm cap_p}\ (fB(x,2\varepsilon),f\overline{B(x,\varepsilon)})
\ge
\left(c_1\frac{d^p(fB(x,\varepsilon))}{[m(fB(x,2\varepsilon))]^{1-n+p}}\right)^{\frac{1}{n-1}}
\end{equation}
where  $c_1$ depends only on $n$ and $p$, $d(A)$ is a diameter and
$m(A)$ is the Lebesgue measure of $A$ in ${\Bbb R}^n.$
\medskip
Combining (\ref{2.4}) and (\ref{2.5}), we obtain that
$$\frac{d(fB(x,\varepsilon))}{\varepsilon}\leq c_2
\left(\frac{m(fB(x,2\varepsilon))}{m(B(x,2\varepsilon))}
\right)^{\frac{1-n+p}{p}}\left(\frac{1}{m(B(x,2\varepsilon))}
\int\limits_{B(x,2\varepsilon)} Q(y)\,
dm(y)\right)^{\frac{n-1}{p}}$$ and hence
$$
L(x,f)\leq\limsup\limits_{\varepsilon\rightarrow
0}\frac{d(fB(x,\varepsilon))}{\varepsilon}\leq c_2\,
\varphi^\frac{1-n+p}{p}(x)Q^{\frac{n-1}{p}}(x)
$$
where
\begin{equation}\label{eq11}
L(x,f)=\limsup\limits_{y\rightarrow x}\frac{|f(y)-f(x)|}{|y-x|}\,.
\end{equation}
Thus,  $L(x,f)<\infty$ a.e. in  $D$. Finally, applying the
Rademacher--Stepanov theorem, see e.g. \cite{Sa}, p. 311, we
conclude that $f$ is differentiable a.e. in $D$.

\begin{corol}{}\label{cor1} Let $D$ be a domain in ${\Bbb R}^n$, $n\ge 2$,
and $f:D\rightarrow \overline{{\Bbb R}^n}$ be a ring
$(p,Q)$--mapping with $Q\in L^1_{loc}$ and $p>n-1$. Suppose that
$f$ is discrete and open. Then the partial derivatives of $f$ are
locally integrable.
\end{corol}

{\it Proof.} Given a compact set $V\subset D,$ we have

$$
\int\limits_{V}L(x,f)\ dx\leq c_2
\int\limits_{V}\varphi^{\frac{1-n+p}{p}}(x)Q^{\frac{n-1}{p}}(x)\
dm(x)$$

Applying the H\"{o}lder inequality, see (17.3) in \cite{BB}, we
obtain

$$
\int\limits_{V}\varphi^{\frac{1-n+p}{p}}(x)Q^{\frac{n-1}{p}}(x)\
dm(x) \leq \left(\int\limits_{V}\varphi(x)\
dm(x)\right)^{\frac{1-n+p}{p}} \left(\int\limits_{V}Q(x)\
dm(x)\right)^{\frac{n-1}{p}}
$$

and since $Q\in L_{loc}^1$
$$
\int\limits_{V}L(x,f)\ dm(x)\leq c_2
N(f,V)^{2/n}\left(\int\limits_{V}Q(x)\
dm(x)\right)^{\frac{n-1}{p}}<\infty\, ,$$ see Lemma 2.3 in
\cite{MRV$_1$}.

\begin{corol}{}\label{cor2}
 Let $D$ be a domain in ${\Bbb R}^n$, $n\ge 2,$
and $f:D\rightarrow \overline{{\Bbb R}^n}$ be a ring
$(p,Q)$--mapping with $Q\in L^1_{loc}$ and $p>n-1$. Suppose that
$f$ is discrete and open. Then
$$\Vert f^{\,\prime}(x)\Vert^p\ \le\ C\cdot|J(x,f)|^{1-n+p}\ Q^{n-1}(x)$$
a.e. where a constant $C$ depends only on $n$ and  $p$.

\end{corol}

\section{On the ACL property of discrete open $(p,Q)$--mappings}\cc
\medskip

\begin{theo}{}\label{th2}
Let $D$ be a domain in ${\Bbb R}^n$, $n\ge 2$, and $f:D\rightarrow
\overline{{\Bbb R}^n}$ be a ring $(p,Q)$--mapping with $Q\in
L^1_{loc}$ and $p>n-1$. Suppose that $f$ is discrete and open.
Then $f\in ACL.$
\end{theo}

{\it Proof.} Without loss of generality we may assume that
$\infty\notin D^{\,\prime}=f(D).$
 Let $I=\{x\in{\Bbb R}^n:a_i<x_i<b_i,\ i=1,\ldots,n\}$
be an  $n$-dimensional interval in $\Rn$ such that
$\overline{I}\subset D$. Then $I=I_0\times J$ where  $I_0$ is an
$(n-1)$-dimensional interval in $\Bbb{R}^{n-1}$ and  $J$ is an open
segment of the axis $x_n$, $J=(a,b)$. Next we identify
$\Bbb{R}^{n-1}\times \Bbb{R}$ with $\Bbb{R}^n$. We prove that for
almost everywhere  segments $J_{z}=\{ z\} \times J\,, \,z\in I_0, $
the mapping $f|_{J_z}$ is absolutely continuous.
\medskip

Consider the set function $\Phi(B)=m\left(f(B\times J)\right)$
defined over the algebra of all Borel sets $B$ in $I_0.$ By 2.2,
2.3 and 2.12 in \cite{MRV$_1$}
\begin{equation}\label{2.1}
\varphi(z)=\limsup\limits_{r\rightarrow
0}\frac{\Phi(B(z,r))}{\Omega_{n-1}r^{n-1}}<\infty
\end{equation}
for a.e. $z\in I_0$ where $B(z,r)$ is a ball in $\R^{n-1}$
centered at the point $z\in I_0$ of the radius $r$ and
$\Omega_{n-1}$ is a volume of the unit ball in ${\Bbb R}^{n-1}$.

Let $\Delta_i, \, i=1,2,...,$ be some  enumeration  $S$ of all
intervals  in $J$ such that  $\overline{\Delta_i}\subset J$ and the
ends of $\Delta_i$ are the rational numbers. Set
$$
\varphi_i(z):=\int\limits_{\Delta_i}Q(z,x_n)\,dx_n.
$$
Then by the Fubini theorem, see e.g. III. 8.1 in \cite{Sa}, the
functions $\varphi_i(z)$  are a.e. finite and integrable in $z\in
I_0$. In addition, by the Lebesgue theorem on differentiability of
the indefinite integral there is  a.e. a finite limit
\begin{equation}\label{2.2}
\lim\limits_{r\rightarrow
0}\frac{\Phi_i(B^{n-1}(z,r))}{\Omega_{n-1}r^{n-1}}=\varphi_i(z)
\end{equation}
where $\Phi_i$ for a fixed $i=1,2,\ldots $ is the set function
$$
\Phi_i(B)=\int\limits_{B}\varphi_i(\zeta)\,d\zeta
$$
given over the algebra of all Borel sets $B$ in $I_0.$
\medskip

Let us show that the mapping  $f$ is absolutely continuous on each
segment $J_z, z\in I_0$, where the finite limits $(\ref{2.1})$ and
$(\ref{2.2})$ exist. Fix one of such a point $z.$ We have to prove
that the sum of diameters of the images of an arbitrary finite
collection of mutually disjoint segments in $J_z=\{z\}\times J$
tends to zero together with the total length of the segments. In
view of the continuity of the mapping $f$, it is sufficient to
verify this fact only for mutually disjoint segments with rational
ends in $J_z$. So, let
$\Delta_i^{*}=\{z\}\times\overline{\Delta_i}\subset J_z$ where
$\Delta_i\in S,\, i=1,...,k$ under the corresponding
re-enumeration of $S$, are mutually disjoint intervals. Without
loss of generality, we may assume that $\overline{\Delta_i}\,
,i=1,...,k$ are also mutually disjoint.
\medskip

Let $\delta>0$ be an arbitrary  rational number which is less than
half of the minimum of the distances between ${\Delta^*_i}\,
,i=1,...,k$, and also less than their distances to the end-points
of the interval $J_z$. Let $\Delta^*_i=\{ z\} \times
[\alpha_i,\beta_i]$ and
$A_i=A_i(r)=B^{n-1}(z,r)\times(\alpha_i-\delta, \beta_i+\delta),$
$i=1,...,k $ where $B^{n-1}(z,r)$ is an open ball in $I_0\subset
\Bbb{R}^{n-1}$  centered at the point  $z$ of the radius $r>0$.

For small $r>0$, $E_i=(A_i,\Delta_i^{*}), i=1,...,k$ are
condensers in $I$ and hence, $fE_i=(fA_i,f\Delta_i^{*}),
i=1,...,k$ are condensers in $D^{\,\prime}.$ By Lemma
\ref{lem2.2},
$${\rm cap_p}\,(fA_i,f\Delta_i^{*})=M_p(\Gamma_{f_{E_i}})\,.$$
Denoting through $\Gamma_{E_i^*}$ a family of maximal $f$--liftings
of $\Gamma_{f_{E_i}}$ starting at $\Delta_i^{*},$ we obtain
$\Gamma_{E_i^*}\subset \Gamma_{E_i}$ and
\begin{equation}\label{eq8}
{\rm cap_p}\,(fA_i,f\Delta_i^{*})\le M_p (f\Gamma_{E_i})\,.
\end{equation}
Let $m$ be a natural number such that $1/m<\delta.$ Consider the
ring $\varepsilon_1<|x-z_0|<\varepsilon_2$ where
$z_0=\left(z,\frac{\alpha_i+\beta_i}{2}\right),$ $\varepsilon_1=
\frac{\beta_i-\alpha_i}{2},$ $\varepsilon_2=
\frac{\beta_i-\alpha_i}{2}+\delta-1/m.$ Let $\Gamma_{i,m}$ is a
path family joining the spheres
$S_1=\left\{|x-z_0|=\varepsilon_1\right\}$ and
$S_2=\left\{|x-z_0|=\varepsilon_2\right\}$ in ${\Bbb R}^n.$ Note
that $\Gamma_{i,m}<\Gamma_{E_i}$ and by (\ref{eq8})
\begin{equation}\label{eq8.1}
{\rm cap_p}\,(fA_i,f\Delta_i^{*})\le M_p (f\Gamma_{i,m})\,.
\end{equation}
Consider the family of functions
$$ \eta_{i,m}(t)\,=\,\left
\{\begin{array}{rr} \frac{1}{r-1/m}, & {\rm если } \ t\in
\left(\frac{\beta_i-\alpha_i}{2},
\frac{\beta_i-\alpha_i}{2}+\delta-1/m\right), \\
0, & {\rm если} t\in {\Bbb
R}\setminus\left(\frac{\beta_i-\alpha_i}{2},
\frac{\beta_i-\alpha_i}{2}+\delta-1/m\right)\,.
\end{array}\right.
$$
as $r<\delta.$ By definition of ring $(p,Q)$--mappings, from
(\ref{eq8.1}) we have
\begin{equation}\label{eq3}
{\rm cap_p}\,(fA_i,f\Delta_i^{*})\leq\,
\frac{1}{(r-1/m)^p}\int\limits_{A_i} Q(x)\, dm(x)\,.
\end{equation}
Letting into the limit in (\ref{eq3}) as $m\rightarrow\infty,$ we
obtain
\begin{equation}\label{3.4}
{\rm cap_p}\,(fA_i,f\Delta_i^{*})\leq\,
\frac{1}{r^p}\int\limits_{A_i} Q(x)\, dm(x)\,.
\end{equation}
On the other hand, by Proposition 6 in \cite{Kr},
\begin{equation}\label{3.5}
{\rm cap_p}\,(fA_i,f\Delta_i^{*})\ge
\left(c\frac{d_i^p}{m_i^{1-n+p}}\right)^{\frac{1}{n-1}}
\end{equation}
where $d_i$ is a diameter of the set $f\Delta_i^{*},$ $m_i$ is a
volume of $fA_i$ and $c$  is a constant depending only on $n$ and
$p$.

Combining $(\ref{3.4})$ and $(\ref{3.5}),$ we have
\begin{equation}\label{3.6}
\left(\frac{d_i^p}{m_i^{1-n+p}}\right)^{\frac{1}{n-1}}\leq
\frac{c_1}{r^p}\int\limits_{A_i}Q(x)dm(x)
\end{equation}
with a constant $c_1$ depending only on $n$, $p$ and all
$i=1,...,k.$

By the discrete H\"{o}lder inequality  see e.g. (17.3) in
\cite{BB}, we obtain
$$\sum\limits_{i=1}^{k} d_i\leq\left(\sum\limits_{i=1}^{k}\left(
\frac{d_i^p}{m_i^{1-n+p}}\right)^{\frac{1}{n-1}}
\right)^{\frac{n-1}{p}}\left(\sum\limits_{i=1}^{k}
m_i\right)^{\frac{1-n+p}{p}}\,,$$
i.e.
$$\left(\sum\limits_{i=1}^{k}
d_i\right)^p\leq\left(\sum\limits_{i=1}^{k}\left(
\frac{d_i^p}{m_i^{1-n+p}}\right)^{\frac{1}{n-1}}
\right)^{n-1}[\Phi(B(z,r))]^{1-n+p}.$$
By $(\ref{3.6})$
$$\left(\sum\limits_{i=1}^{k} d_i\right)^p\leq c_2
\left[\frac{\Phi(B(z,r))}{\Omega_{n-1}r^{n-1}}\right]^{1-n+p}\left(\sum\limits_{i=1}^{k}\frac{\int\limits_{A_i}Q(x)\,dm(x)}
{\Omega_{n-1}r^{n-1}}\right)^{n-1}$$
where $c_2$ depends only on $n$ and $p$. Passing to the limit
first as $r\rightarrow 0$ and then as $\delta\rightarrow 0$, we
obtain
\begin{equation}\label{3.10}
\left(\sum\limits_{i=1}^{k} d_i\right)^p\leq c_2\,
[\varphi(z)]^{1-n+p}\left(\sum\limits_{i=1}^{k}\varphi_i(z)
\right)^{n-1}\,.\end{equation} Finally, in view of (\ref{3.10}),
the absolute continuity of the indefinite integral of $Q$ over the
segment  $J_z$ implies the  absolute continuity of the mapping $f$
over the same segment. Hence $f\in ACL$.

\medskip

Combining Theorem \ref{th2} and Corollary \ref{cor1}, we obtain a
following conclusion, see also \cite{Ma}.

\begin{corol}{}\label{cor3} Let $D$ be a domain in ${\Bbb R}^n$, $n\ge 2$, and $f:D\rightarrow
\overline{{\Bbb R}^n}$ be a ring $(p,Q)$--mapping with $Q\in
L^1_{loc}$ and $p>n-1$. Suppose that $f$ is discrete and open.
Then $f\in W^{1,1}_{loc}$.
\end{corol}

\medskip
\noindent Ruslan Salimov and Evgeny Sevost'yanov, \\
Institute of Applied Mathematics and Mechanics,\\
National Academy of Sciences of Ukraine, \\
74 Roze Luxemburg str., 83114 Donetsk, UKRAINE \\
Phone: +38 -- (062) -- 3110145 Fax: \ \ \ \,+38 -- (062) -- 3110285 \\
salimov@iamm.ac.donetsk.ua, e\_sevostyanov@rambler.ru

\end{document}